\documentclass[12pt,reqno]{amsart}

\usepackage{amsfonts, amssymb, curves, epic, epsf, amsmath}
\newcommand{\R}{\mathbb{R}}

\newcommand{\eps}{\epsilon}

\newcommand{\larr}{\left( \begin{array}{c}}
\newcommand{\rarr}{\end{array} \right) }

\newcommand{\lsqarr}{\left[ \begin{array}{c}} 
\newcommand{\rsqarr}{\end{array} \right]}

\begin{document}

\newtheorem{theorem}{Theorem}
\newtheorem*{defn}{Definition}
\newtheorem{cor}[theorem]{Corollary}
\newtheorem{corollary}[theorem]{Corollary}
\newtheorem{lemma}[theorem]{Lemma} \newtheorem{prop}[theorem]{Proposition}
\newtheorem{example}[theorem]{Example}

\title[An approximation theorem for tiling maps]{An approximation 
theorem for maps between tiling spaces}
\author{Betseygail Rand}\address{Department of Mathematics, 
Texas Lutheran University, Seguin, TX 78155}\email{brand@tlu.edu}
\author{Lorenzo Sadun}\address{Department of Mathematics, 
University of Texas, Austin, TX 78712} \email{sadun@math.utexas.edu}

 \begin{abstract} 
We show that every continuous map from one translationally finite tiling
space to another can be approximated by a local map. If two local maps 
are homotopic, then the homotopy can be chosen so that every interpolating
map is also local.
\end{abstract} 

 \maketitle

\section{Introduction}

Many aspects of tiling theory, such as pattern-equivariant cohomology
\cite{K, KP}, are built around local data. It might not matter what a
tiling looks like near infinity, but it might matter crucially that a
certain tile sits exactly {\em here}. Unfortunately, maps between tiling spaces
may not preserve local data. Even topological conjugacies need not be local
maps \cite{Pet,RS}. To show that structures built from local data are actually
topological invariants, we need to show that arbitrary maps between tiling
spaces can be approximated by local maps, and that homotopies between local
maps can be chosen to preserve locality at all times. These are the two main
theorems of this paper. In addition, we show that all constructions can be
chosen to preserve whatever discrete rotational symmetry exists. 

For our purposes, a {\em tiling} is a decomposition of the plane (or, more
generally, of $\R^d$) into a countable union of closed polygons (or
polyhedra) that overlap only on their boundaries.  These polygons are
called {\em tiles}. One can consider
more complicated shapes than polygons, but there is a standard trick,
involving Voronoi cells \cite{Priebe}, that converts non-polygonal
tilings into polygonal tilings with the same mathematical properties.

A {\em patch} is a sub-collection of tiles in the tiling. For any
tiling $T$ and set $S \subset \R^d$, we let $[S]^T$ be the set of all
tiles in $T$ that intersect $S$. The {\em central patch of radius $R$} is
the patch defined by $S=B_R(0)$, the closed 
ball of radius $R$ around the origin. 
For any $x\in \R^d$, $T-x$ is a translate of $T$; a neighborhood of the 
origin in $T-x$ looks like a neighborhood of the point $x$ in $T$. 

Two tilings are considered $\epsilon$-close if they agree on $B_{1/\eps}(0)$,
up to a rigid motion that moves points in $B_{1/\eps}(0)$ by $\epsilon$ or
less. This rigid motion need not be a translation, but in most
examples it is. A {\em tiling space} is a set of tilings that is
complete in the tiling metric and is invariant under translations. The
closure of the set of translates of any given tiling $T$ is special
kind of tiling space, called the {\em hull} of $T$. 

A tiling space has {\em translational finite local complexity}, or is
{\em translationally finite}, if the set of all patches of radius $R$
is finite up to translation.  A {\em tiling} is translationally finite
if its hull is translationally finite.  Translationally finite tiling
spaces are necessarily compact.

\begin{defn}
If $\Omega_1$ and $\Omega_2$ are tiling spaces, we say a map 
$f: \Omega_1 \to \Omega_2$ is {\em local with radius $R$} if, whenever
two tilings $T, T' \in \Omega_1$ have identical central patches of radius $R$,
then $f(T)$ and $f(T')$ have identical central patches of radius 1. 
\end{defn}
In other
words, if $f$ is local, 
you don't need to know the behavior of $T$ near infinity to specify
the behavior of $f(T)$ near the origin.  

For factor maps (i.e., maps that
commute with translation), this is the analog of a sliding block code. 
On subshifts, continuous factor maps are always sliding block
codes \cite{LiM}, but on tiling spaces, continuous 
factor maps need not be local \cite{Pet, RS}. The problem has
to do with small rigid motions. If $T$ and $T'$ agree on a large ball, then
$T$ and $T'$ are close in the tiling metric, which means that $f(T)$ and 
$f(T')$ are close, which means that $f(T)$ and $f(T')$ agree on a large ball,
{\em up to a small motion}. 

In this paper, we show how to get rid of the small motion, although 
typically at the cost of not remaining a factor map. Indeed, we do 
not assume that our maps commute
with translation to begin with! We merely show how to approximate
arbitrary continuous maps with local continuous maps (Section 2), and
how to approximate arbitrary homotopies between local maps with
homotopies that preserve locality.  These results are extensions, with
streamlined proofs, of results first announced in \cite{Rand}. 

Note that these results apply only to {\em translationally} finite tilings. 
If a tiling has finitely many patches of radius $R$ up to Euclidean motion,
but not up to translation (e.g., the pinwheel tiling \cite{pin}), then the
averaging trick used to prove Theorem 1 breaks down, since we would have to 
average elements of a non-Abelian group. For a discussion of what can be 
proved for tilings with (rotational) finite local complexity, see \cite{Rand}.

We thank Franz Gaehler, John Hunton, Johannes Kellendonk, and Ian
Putnam for useful discussions.  The work of the second author is
partially supported by the National Science Foundation.

\section{The approximation theorem} 

\begin{theorem}
Let $\Omega_1$ be the hull of a translationally finite
and non-periodic tiling $T_0$, let $\Omega_2$ be a translationally finite
tiling space,
and let $f: \Omega_1 \to \Omega_2$ be a continuous map. 
For each $\eps > 0$
there exists a continuous local map $f_\eps: \Omega_1 \to \Omega_2$ such
that $f$ and $f_\eps$ differ only by a small translation. Specifically, there
exists a continuous function $s_\eps: \Omega_1 \to \R^d$ such that,
for each tiling $T$, $f_\eps(T) = f(T) - s_\eps(T)$ and $|s_\eps(T)|<\eps$.
\end{theorem}

\begin{proof}
$\Omega_1$ is the hull of $T_0$ and is
translationally finite. This means that each patch $P$ of a tiling
$T \in \Omega_1$ is found somewhere in $T_0$, say at position $x$, 
so $P$ is the central patch of $T_0-x$.
Since $f$ is continuous and $\Omega_1$ is compact, $f$ is uniformly 
continuous. Pick $\delta$ such that, if two tilings 
$T_1$, $T_2$ agree on $B_{1/\delta}(0)$, then $f(T_1)$ and $f(T_2)$ are within
$\eps$. Let $R=2/\delta$. We will construct $f_\eps$ to be local with 
radius $R+\delta$.

$T_0$ is a single point in $\Omega_1$, but it is sometimes convenient to view
$T_0$ as a marked copy of $\R^d$.
For $x,y \in \R^d$, let $x \sim y$ if $[B_R(0)]^{T_0-x} =
[B_R(0)]^{T_0-y}$. In other words, $x \sim y$ if the patch of radius
$R$ around $x$ in $T_0$ looks like the patch around $y$. Let $K_R$ be
the quotient of $\R^d$ by this equivalence relation. $K_R$ is a
branched $d$-manifold \cite{BDHS,S} that parametrizes the possible
patches of radius $R$. Since every patch of every tiling is found
somewhere in $T_0$, there is a natural projection $\pi: \Omega_1 \to
K_R$ that sends each tiling to the description of its central patch.

$K_R$ is a CW complex \cite{BDHS}, 
and is easily decomposed into disjoint cells of
dimension up to $d$. For each cell $C$, pick a connected region
$\tilde C \subset \R^d$ that represents $C$. That is, each point $p \in C$
is the equivalence class of a unique point $h(p) \in \tilde C$. 
Restricted to a single cell $C$, the map $h: K_R \to \R^d$ is continuous, 
but $h$ may jump as we pass from one cell to another. Let $g = h \circ \pi:
\Omega_1 \to \R^d$. 

Let $T$ be any tiling in $\Omega_1$. Since $\pi(T)$ describes the central 
patch of $T$, and since $g(T)$ is a point in $T_0$ whose patch of
radius $R$ agrees with the central patch of $T$, $T$ and $T_0 - g(T)$
agree exactly on $B_R(0)$.  This implies that $f(T)$ agrees with 
$f(T_0-g(T))$ on $B_{1/\eps}(0)$, up to 
translation by up to $\eps$. If $\eps$ is small, this
translation is unique. Let $\tilde s_\eps(T)$ be the unique small element
of $\R^d$ such that $f(T) - \tilde s_\eps(T)$ agrees exactly with
$f(T_0-g(T))$ on $B_R(0)$. Let $\tilde f_\eps(T) = f(T) - \tilde s_\eps(T)$. 

By construction, $\tilde f_\eps$ is local, but it may not be
continuous, since $h$ may have jump discontinuities. We remedy this by
convolving $\tilde f_\eps$ with a bump function, insofar as the
convolution of a smooth function with a step function is smooth.  Let
$\phi: \R^d \to \R$ be a smooth function of total integral 1,
supported on $B_\delta(0)$.  For each $y \in B_\delta(0)$, $\tilde
f_\eps(T-y)$ is a small translate of $f(T-y)$, and hence a small
translate of $f(T)$. Let $\rho(T,y)$ be the unique small element of
$\R^d$ such that $\tilde f_\eps(T-y) = f(T) - \rho(T,y)$. Let
\begin{equation}
s_\eps(T) = \int \phi(y) \rho(T,y) dy, \qquad \hbox{and} 
\end{equation}
\begin{equation} 
f_\eps(T) = f(T)- s_\eps(T).
\end{equation}

It is clear that $f_\eps$ is continuous along a translational orbit. 
What remains is to show that $f_\eps$ is local. Suppose that
$T_1$ and $T_2$ agree on $B_{R+\delta}$.  $T_1-y$ and $T_2-y$
agree on $B_R(0)$, so $\tilde f_\eps(T_1-y)$ and $\tilde f_\eps(T_2-y)$
agree on a central patch. This means that $\rho(T_1,y) - \rho(T_2,y)=
\alpha$, where $\alpha$ is 
the translation needed to take the central patch of $f(T_1)$
onto the central patch of $f(T_2)$.  Integrating over $y$, we obtain
$s_\eps(T_1) - s_\eps(T_2) = \alpha$, so the central
patch of $f_\eps(T_1)$ agrees exactly with the central patch of 
$f_\eps(T_2)$. 
\end{proof}

\section{The homotopy theorem}

\begin{theorem}
Let $\Omega_{1,2}$ be as before. 
Let $f_0: \Omega_1 \to \Omega_2$ and $f_1: \Omega_1 \to \Omega_2$ be 
local maps. If $F: [0,1] \times \Omega_1 \to \Omega$ is a homotopy
between $f_0$ and $f_1$, then there is another homotopy $\tilde F$
between $f_0$ and $f_1$ such that each time slice is a local map
from $\Omega_1$ to $\Omega_2$. 
\end{theorem}

\begin{proof}
We apply the method outlined in the proof of Theorem 1 to each time slice 
$f_t$. Since the unit interval is compact, for any $\eps$ one can choose 
values of $\delta$ and $R$ that work for every $f_t$. The resulting family
of local maps $\tilde F: [0,1] \times \Omega_1 \to \Omega$ gives a homotopy
between $f_{0,\eps}$ and $f_{1,\eps}$. What remains is to construct a (local)
homotopy between $f_0$ and $f_{0,\eps}$, and likewise between $f_1$ and 
$f_{1,\eps}$. 

Since $f_0$ is already local, $\tilde f_{0,\eps} = f_0$, and $\rho(T,0)=0$
for every tiling $T$.  For each $t>0$, let $\phi_t(y)=t^{-d}\phi(y/t)$.
When $t=1$, we have our usual function $\phi$, and equation (1)
gives $f_{0,\eps}(T)$. 
As $t \to 0$, $\phi_t$ becomes a delta function, the integral approaches
zero, and equation (1) gives a limiting value of $f_0(T)$. The same argument
gives a local homotopy between $f_1$ and $f_{1,\eps}$.
 
\end{proof}

\section{Rotations}

Up to now we have been discussing tilings and the action of the translation
group on them. For many tilings, such as the Penrose tiling, the rotation
properties are also interesting.  The following theorem extends
Theorems 1 and 2 to that setting.

\begin{theorem} Let $\Omega_1$ and $\Omega_2$ be hulls of 
translationally finite tilings. Suppose that a finite
subgroup $G$ of $SO(d)$ acts naturally on $\Omega_1$ and $\Omega_2$, and 
suppose that $f: \Omega_1 \to \Omega_2$ is a continuous map that
intertwines the actions of $G$. Then the approximation
$f_\eps$ of Theorem 1, besides being local, can be chosen 
to interwine the action of $G$. If $f_0$ and $f_1$ are homotopic maps
$\Omega_1 \to \Omega_2$, and if each is local and each
intertwines the action of $G$, then the
homotopy between them can be chosen so that each $f_t$ is local and 
intertwines the action of $G$. 
\end{theorem}

\begin{proof}
  Construct $s_\eps(T)$ exactly as in Theorem 1, only with a
  rotationally symmetric function $\phi(y)$, and then average over the
  group, defining $\bar s_\eps(T) = \frac{1}{|G|} \sum_{g \in G}
  g^{-1} s_\eps(g T)$ and $f_\eps(T) = f(T) - \bar
  s_\eps(T)$. That proves the first half of the theorem. 
Applying the same construction to the homotopy between $f_0$ and $f_1$ proves
the second half of the theorem. 
\end{proof}

\end{document}